\documentclass[11pt]{amsart}
\usepackage{amssymb}
\usepackage{url}
\usepackage{enumerate}
\usepackage{tikz}
\newtheorem{lemma}{\bf Lemma}
\newtheorem{theorem}[lemma]{\bf Theorem}

\newtheorem{proposition}[lemma]{\bf Proposition}

\DeclareMathOperator{\h}{h}

\begin{document}
\parskip = 0mm
\title[An Issue Raised by Duffus in 1978]{An Issue Raised in 1978 by a Then-Future Editor-in-Chief of the Journal {\sl Order}\\ \small Does the Endomorphism Poset of a Finite Connected Poset Tell Us That the Poset Is Connected?}
\author{Jonathan David Farley}
\address{Department of Mathematics, Morgan State University, 1700 E. Cold Spring Lane, Baltimore, MD 21251, United States of America, {\tt lattice.theory@gmail.com}}



\begin{abstract}
In 1978, Dwight Duffus---editor-in-chief of the journal {\sl Order} from 2010 to 2018 and chair of the Mathematics Department at Emory University from 1991 to 2005---wrote that ``it is not obvious that $P$ is connected and $P^P$ isomorphic to $Q^Q$ implies that $Q$ is connected,'' where $P$ and $Q$ are finite non-empty posets.  We show that, indeed, under these hypotheses $Q$ is connected and $P\cong Q$.  
\end{abstract}

\thanks {The author would like to thank Bernd Schr\"oder for his detailed critique of the paper.  The author would like to thank Dwight Duffus for sending him a reference and sections of his thesis (which the author has had since graduate school but did not have handy when writing this paper) and Dominic van der Zypen and Zeinab Bandpey for their willingness to look at these results.}

\maketitle


\def\Qa{\mathbb{Q}_0}
\def\Qb{\mathbb{Q}_1}
\def\Q{\mathbb{Q}}
\def\card{{\rm card}}
\parskip = 2mm
\parindent = 10mm
\def\Part{{\rm Part}}
\def\P{{\mathcal P}}
\def\Eq{{\rm Eq}}
\def\cld{Cl_\tau(\Delta)}
\def\Csing{{\mathcal C}_{\{*\}}}
\def\Cftwo{{\mathcal C}_{{\rm fin}\rangle1}}
\def\Cinf{{\mathcal C}_{\infty}}
\def\Pcf{{\mathcal P}_{\rm cf}}
\def\Fn{{\mathcal F}_n}
\def\proof{{\it Proof. }}


\vspace*{-4mm} 

In the 1979 {\sl Proceedings of the American Mathematical Society}, Duffus and Wille proved that $P^P\cong Q^Q$ implies $P\cong Q$ if $P$ and $Q$ are both finite, non-empty, and connected \cite[Theorem, p.14]{DufWilGI}.  Duffus---editor-in-chief of the journal {\sl Order} from 2010 to 2018 and chair of the Mathematics Department at Emory University from 1991 to 2005---said in 1984 \cite[p. 90]{DufHDa}, ``It is still an open problem to show connectedness can be dropped.'' He notes in his 1978 thesis \cite[p. 53, as told to the author by Duffus]{DufGH}, ``[I]t is not obvious that $P$ is connected and $P^P$ isomorphic to $Q^Q$ implies that $Q$ is connected.'' (Note that, if this is the case, then the Duffus--Wille result implies $P\cong Q$.)  We prove that, indeed, if $P$ and $Q$ are finite, non-empty posets such that $P^P\cong Q^Q$ and $P$ is connected, then $Q$ is connected.

In other words, we have resolved the issue from Duffus's 1978 thesis in that we have ``half-dropped'' the connectedness hypothesis used in the 1979 {\sl Proceedings of the American Mathematical Society} paper, whereas the 1984 problem asked if it could be dropped entirely.

\medskip

We assume the reader is familiar with the basic facts about the arithmetic of ordered sets, the basic consequences of Hashimoto's Refinement Theorem, and Professor Birkhoff's theorem on finite distributive lattices (e.g., \cite[Propositions 3.1 and 4.1]{McKJC} and \cite[Theorems 5.9 and 5.12]{DavPriJB}).  A good reference is \cite{SchJC}.

Let $P$ and $Q$ be posets.  For $p,p'\in P$, we write $p\equiv p'$ if $p$ and $p'$ are in the same connected component of $P$; $\h(P)$ denotes the height of the finite, non-empty poset $P$, the largest value of $|C|-1$ for $C$ a chain in $P$, and $\h_P(p)$ denotes the height of an element $p$ in $P$. Note that $\h(P\times Q)=\h(P)+\h(Q)$ for $P$ and $Q$ finite and non-empty \cite[Chapter I, \S9, Exercise 4(a)]{BirDH}.

Let $P^Q$ denote the poset of order-preserving maps from $Q$ to $P$, where, for $f,g\in P^Q$,
$$
f\le g\ \text{if for all}\ q\in Q\text{,}\ f(q)\le_Q g(q)
$$
\noindent\cite[p. 312]{BirCG}.  For $p\in P$, denote the constant map $f(q)=p$ for all $q\in Q$ by $\langle p\rangle$.  Let $\mathcal D(P^Q)$ denote
$$
\{g\in P^Q\mid g\ \text{is constant on each connected component of}\ Q\}
$$
and for $Q\ne\emptyset$ let $\mathcal C(P^Q)$ denote
$$
\{f\in P^Q\mid f\equiv g\ \text{for some}\ g\in\mathcal D(P^Q)\}
$$
\noindent\cite[Definition 2.1]{McKJC}.  Note that, despite the notation, $\mathcal C(P^Q)$ is really a function of $P$ and $Q$, not $P^Q$ (as far as the author knows).

\begin{center}

    \begin{tikzpicture}[scale=.85]

    \draw[fill] (-4,0) circle (.05cm);
    \draw[fill] (4,0) circle (.05cm);
    \draw[fill] (-7,1) circle (.05cm);
    \draw[fill] (-5,1) circle (.05cm);
    \draw[fill] (-3,1) circle (.05cm);
    \draw[fill] (-1,1) circle (.05cm);
    \draw[fill] (7,1) circle (.05cm);
    \draw[fill] (5,1) circle (.05cm);
    \draw[fill] (3,1) circle (.05cm);
    \draw[fill] (1,1) circle (.05cm);

    \draw[fill] (-7,2) circle (.05cm);
    \draw[fill] (-6,2) circle (.05cm);
    \draw[fill] (-5,2) circle (.05cm);
    \draw[fill] (-3,2) circle (.05cm);
    \draw[fill] (-2,2) circle (.05cm);
    \draw[fill] (-1,2) circle (.05cm);

    \draw[fill] (7,2) circle (.05cm);
    \draw[fill] (6,2) circle (.05cm);
    \draw[fill] (5,2) circle (.05cm);
    \draw[fill] (3,2) circle (.05cm);
    \draw[fill] (2,2) circle (.05cm);
    \draw[fill] (1,2) circle (.05cm);

    \draw[fill] (-7,3) circle (.05cm);
    \draw[fill] (-5,3) circle (.05cm);
    \draw[fill] (-3,3) circle (.05cm);
    \draw[fill] (-1,3) circle (.05cm);

    \draw[fill] (7,3) circle (.05cm);
    \draw[fill] (5,3) circle (.05cm);
    \draw[fill] (3,3) circle (.05cm);
    \draw[fill] (1,3) circle (.05cm);

    \draw[fill] (-4,4) circle (.05cm);
    \draw[fill] (4,4) circle (.05cm);

    \draw[fill] (8,2) circle (.05cm);
    \draw[fill] (9,2) circle (.05cm);
    \draw[fill] (10,2) circle (.05cm);
    \draw[fill] (11,2) circle (.05cm);

    \draw (-4,0) -- (-7,1);
    \draw (-4,0) -- (-5,1);
    \draw (-4,0) -- (-3,1);
    \draw (-4,0) -- (-1,1);

    \draw (4,0) -- (7,1);
    \draw (4,0) -- (5,1);
    \draw (4,0) -- (3,1);
    \draw (4,0) -- (1,1);

    \draw (-7,1) -- (-7,2);
    \draw (-7,1) -- (-2,2);

    \draw (-5,1) -- (-7,2);
    \draw (-5,1) -- (-3,2);

    \draw (-3,1) -- (-3,2);
    \draw (-3,1) -- (1,2);

    \draw (-1,1) -- (-2,2);
    \draw (-1,1) -- (1,2);

    \draw (7,1) -- (7,2);
    \draw (7,1) -- (6,2);

    \draw (5,1) -- (7,2);
    \draw (5,1) -- (3,2);

    \draw (3,1) -- (3,2);
    \draw (3,1) -- (-1,2);

    \draw (1,1) -- (6,2);
    \draw (1,1) -- (-1,2);

    \draw (-7,2) -- (-7,3);
    \draw (-7,2) -- (-5,3);

    \draw (-6,2) -- (-7,3);
    \draw (-6,2) -- (-1,3);

    \draw (-5,2) -- (-5,3);
    \draw (-5,2) -- (-3,3);

    \draw (-1,2) -- (-3,3);
    \draw (-1,2) -- (-1,3);

    \draw (1,2) -- (1,3);
    \draw (1,2) -- (3,3);

    \draw (2,2) -- (1,3);
    \draw (2,2) -- (7,3);

    \draw (5,2) -- (3,3);
    \draw (5,2) -- (5,3);

    \draw (7,2) -- (5,3);
    \draw (7,2) -- (7,3);

    \draw (-7,3) -- (-4,4);
    \draw (-5,3) -- (-4,4);
    \draw (-3,3) -- (-4,4);
    \draw (-1,3) -- (-4,4);

    \draw (7,3) -- (4,4);
    \draw (5,3) -- (4,4);
    \draw (3,3) -- (4,4);
    \draw (1,3) -- (4,4);

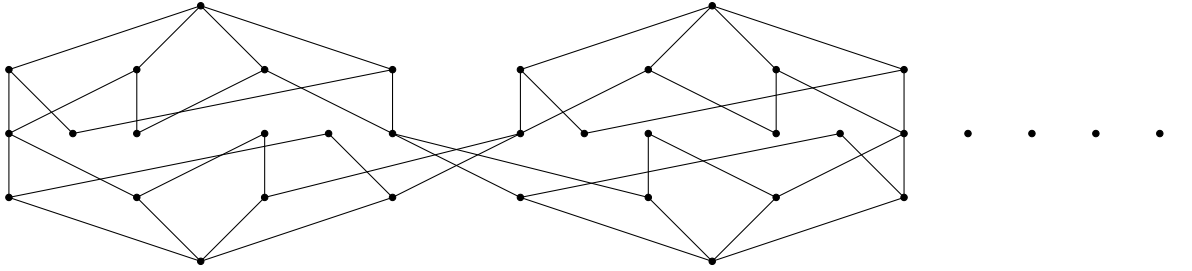
    \draw (0,-1) node {\bf Figure. $P^P$ where $P$ is the 4-element crown \cite{DufGH}};

    \end{tikzpicture}
    
\end{center}

The $n$-element chain is denoted $\bf n$.  A poset $P$ is {\it directly irreducible} if $|P|\ne1$ and whenever $P\cong A\times B$ for posets $A$ and $B$, then $|A|=1$ or $|B|=1$.  A finite poset $P$ is {\it absolutely $\mathcal C$-indecomposable} if it is connected, directly irreducible, and, whenever $P\cong\mathcal C(A^B)$ for posets $A$ and $B$, $B\ne\emptyset$, we have $P\cong A$ and $|B|=1$.

The following result comes from a remark McKenzie leaves to the reader after \cite[Proposition 4.1]{McKJC}, which nonetheless is true for finite, non-empty posets.  (For a proof of (1), see \cite[Theorem 2.8]{FarBJa}.)

\begin{proposition} Let $P$, $Q$, and $R$ be non-empty posets.  Then:
\begin{itemize}
\item[(1)] $\mathcal C(P^{Q\times R})\cong\mathcal C\big(\mathcal C(P^Q)^R\big)$ if $P$, $Q$, and $R$ are finite;
\item[(2)] $\mathcal C(P^{Q+R})\cong\mathcal C(P^Q)\times\mathcal C(P^R)$; 
\item[(3)] if $R$ is connected, $\mathcal C\big((P+Q)^R\big)\cong\mathcal C(P^R)+\mathcal C(Q^R)$;
\item[(4)] $\mathcal C[(Q\times R)^P]\cong\mathcal C(Q^P)\times\mathcal C(R^P)$. $\qed$
\end{itemize}
\end{proposition}

We will also use the colossal structure theorems of McKenzie \cite[Theorems 8.1, 9.1, and 9.2]{McKJC}:

\begin{theorem} \begin{itemize} \item[(1)] Let $A$, $B$, $C$, and $D$ be finite, non-empty, connected posets.  Assume $C$ and $D$ are directly irreducible and $C\ncong D$.  Then if $\mathcal C(A^C)\cong\mathcal C(B^D)$, there exists a finite, non-empty, connected poset $E$ such that $A\cong\mathcal C(E^D)$ and $B\cong\mathcal C(E^C)$. (Note that ``$\bf A^Q$'' should be ``$\bf E^Q$'' on \cite[p. 211]{McKJC}.)
\item[(2)] Let $A$, $B$, and $C$ be finite, non-empty posets.  Assume $C$ is connected.  Then $A\cong B$ if $\mathcal C(A^C)\cong\mathcal C(B^C)$.
\item[(3)] Let $n\in\mathbb N$. Let $A$ and $B$ be finite, non-empty, connected posets, and let $C_1,\dots,C_n$ be posets such that $\mathcal C(A^B)\cong C_1\times\cdots\times C_n$. Then there exist finite, non-empty posets $A_1,\dots,A_n$ such that $C_i\cong\mathcal C(A_i^B)$ for $i=1,\dots,n$ and $A\cong A_1\times\cdots\times A_n$.$\qed$
\end{itemize}
\end{theorem}

We now extend Theorem 2 to the analogue of \cite[Theorem 8.2]{JonMcKHB}, using similar steps.

\begin{lemma} Let $A$, $B$, $C$, and $D$ be finite, non-empty, connected posets such that $\mathcal C(A^C)\cong\mathcal C(B^D)$.  Assume that no non-trivial poset is isomorphic to both a direct factor of $C$ and a direct factor of $D$.  Then there exists a finite, non-empty, connected poset $E$ such that $A\cong\mathcal C(E^D)$ and $B\cong\mathcal C(E^C)$.
\end{lemma}

\proof We use induction. Let $C=C_1\times\cdots\times C_n$ and $D=D_1\times\cdots\times D_m$ be products of directly irreducible posets, where $m,n\ge0$.  Since $C$ and $D$ are connected, so are $C_i$ and $D_j$ ($i=1,\dots,n$; $j=1,\dots,m$).

{\it Case 1. $m=0$ or $n=0$}

If $n=0$, then $\mathcal C(A^C)\cong A$, so let $E=B$.

{\it Case 2. $m=1=n$}

Use Theorem 2.

{\it Case 3. $m+n\ge3$ and $m,n\ge1$}

Without loss of generality, $n\ge2$.  As $\mathcal C(A^{C_1\cdot\cdots\cdot C_n})\cong\mathcal C(B^D)$, by Proposition 1(1) $\mathcal C\big(\mathcal C(A^{C_1})^{C_2\cdot\cdots\cdot C_n}\big)\cong\mathcal C(B^D)$.  Since $A,C\ne\emptyset$, we have $A^{C_1}\ne\emptyset$ and $\mathcal C(A^{C_1})\ne\emptyset$.  Also, $\mathcal C(A^{C_1})$ is connected since $A\ne\emptyset$ is connected.  Thus there exists a finite, non-empty, connected poset $E_1$ such that $\mathcal C(A^{C_1})\cong\mathcal C(E_1^D)$ and $B\cong\mathcal C(E_1^{C_2\cdot\cdots\cdot C_n})$.  Hence there exists a finite, non-empty, connected poset $E$ such that $A\cong\mathcal C(E^D)$ and $E_1\cong\mathcal C(E^{C_1})$.  Therefore by Proposition 1(1), $B\cong\mathcal C\big(\mathcal C(E^{C_1})^{C_2\cdot\cdots\cdot C_n}\big)\cong\mathcal C(E^{C_1\cdot\cdots\cdot C_n})\cong\mathcal C(E^C)$.$\qed$

\begin{theorem} Let $A$, $B$, $C$, and $D$ be finite, non-empty, connected posets such that $\mathcal C(A^C)\cong\mathcal C(B^D)$.  Then there exist finite, non-empty, connected posets $E$, $X$, $Y$, and $Z$ such that $A\cong\mathcal C(E^X)$, $B\cong\mathcal C(E^Y)$, $C\cong Y\times Z$, and $D\cong X\times Z$.
\end{theorem}

\proof Hashimoto's Refinement Theorem tells us we can find finite, non-empty, connected posets $X$, $Y$, and $Z$ such that $C\cong Y\times Z$, $D\cong X\times Z$, and $X$ and $Y$ do not share a non-trivial direct factor.  Thus, by Proposition 1(1),
$$
\mathcal C\big(\mathcal C(A^Y)^Z\big)\cong\mathcal C(A^{Y\times Z})\cong\mathcal C(A^C)\cong\mathcal C(B^D)\cong\mathcal C(B^{X\times Z})\cong\mathcal C\big(\mathcal C(B^X)^Z\big).
$$
\noindent Now $\mathcal C(A^Y),\mathcal C(B^X)\ne\emptyset$ since $A$, $B$, $X$, and $Y$ are non-empty and connected.  Hence, by Theorem 2(2), $\mathcal C(A^Y)\cong\mathcal C(B^X)$.  By Lemma 3, there exists a finite, non-empty, connected poset $E$ such that $A\cong\mathcal C(E^X)$ and $B\cong\mathcal C(E^Y)$.$\qed$

In \cite[Theorem 3.2.6]{DufGH}, Duffus proves an outstanding version of Hashimoto's Refinement Theorem for posets $A$, $B$, $C$, and $D$ that are sums of connected posets with a finite maximal chain such that $A\times B\cong C\times D$ and $A$ and $C$ are connected. He never published this proof, although it is used in \cite{DufWilGI}, so we provide a proof of an $\epsilon$-extension of his refinement theorem now:  We only require one of the four posets to be connected. (We don't actually need this result, but the structure of its proof will help the reader understand the proof of the theorem we do need: this proof foreshadows the proof of Theorem 8.)

\begin{theorem} Let $A$, $B$, $C$, and $D$ be finite posets. Let $A$ be connected and non-empty.  Assume that $A\times B\cong C\times D$.  Then there exist posets $W$, $X$, $Y$, and $Z$ such that $A\cong W\times X$, $B\cong Y\times Z$, $C\cong W\times Y$, and $D\cong X\times Z$.
\end{theorem}

\proof Note that $B=\emptyset$ if and only if either $C=\emptyset$ or $D=\emptyset$.  First assume $B=\emptyset=C$. Then let $W=A$, $X={\bf 1}$, $Y=\emptyset$, and $Z=D$.  Next, assume $B=\emptyset=D$.  Then let $W={\bf 1}$, $X=A$, $Y=C$, and $Z=\emptyset$.

From now on, assume that $B$, $C$, and $D$ are non-empty.

Let $B=\sum_{b\in H} B_h$, $C=\sum_{i\in I} C_i$, and $D=\sum_{j\in J} D_j$ be the decompositions of $B$, $C$, and $D$, respectively, into connected components.  Since $A\times B\cong \sum_{h\in H} A\times B_h$ has the same number of components as $C\times D\cong\sum_{\substack{i\in I\\j\in J}} C_i\times D_j$, there is a bijection $\Psi:I\times J\to H$ such that $C_i\times D_j\cong A\times B_{\Psi(i,j)}$ for all $i\in I$, $j\in J$.  

Let $A_1,A_2,\dots,A_r$ ($r\ge0$) be the pairwise non-isomorphic connected, directly indecomposable posets that could arise in the factorizations of any of $A$, $B_h$ ($h\in H$), $C_i$ ($i\in I$), and $D_j$ ($j\in J$).

Say $A\cong\prod_{\ell=1}^r(A_\ell)^{k_\ell}$ where $k_\ell\ge0$ ($\ell=1,\dots,r$).

For each $\ell\in\{1,\dots,r\}$, let $c_\ell\in\mathbb N_0$ be the highest power of $A_\ell$ such that $A_\ell^{c_\ell}$ is a factor of all of the $C_i$ ($i\in I$).  Let $W:=\prod_{\ell=1}^r A_\ell^{\min\{c_\ell,k_\ell\}}$.  Let $X:=\prod_{\ell=1}^r A_\ell^{k_\ell-\min\{c_\ell,k_\ell\}}$.  Then $W\times X\cong\prod_{\ell=1}^r A_\ell^{k_\ell}\cong A$.

Clearly, for all $i\in I$, $W$ is a factor of $C_i$.

Hence, we may let $\widetilde{C_i}$ be $C_i$ with $W$ factored out---that is, $C_i\cong\widetilde{C_i}\times W$ ($i\in I$).  Let $Y:=\sum_{i\in I}\widetilde{C_i}$.

\quad {\it Claim 1. For each $j\in J$, $X$ is a factor of $D_j$.}

{\it Proof of claim.} Assume for a contradiction that $X$ does not divide $D_j$ for some $j\in J$.  Then there exists $\ell\in\{1,\dots,r\}$ such that $A_\ell^{k_\ell-c_\ell}$ does not divide $D_j$ and so $k_\ell>c_\ell$.  Pick $i\in I$ such that $c_\ell$ is the highest power of $A_\ell$ dividing $C_i$.  Then the highest power of $A_\ell$ dividing $C_i\times D_j$ is less than $c_\ell+k_\ell-c_\ell=k_\ell$, a contradiction, since $C_i\times D_j\cong A\times B_{\Psi(i,j)}$ and $A_{\ell}^{k_\ell}$ divides the right-hand side. $\qed$

By Claim 1, we may let $\widetilde{D_j}$ be $D_j$ with $X$ factored out ($j\in J$).  Let $Z=\sum_{j\in J} \widetilde{D_j}$.

\quad {\it Claim 2. For all $(i,j)\in I\times J$, $B_{\Psi(i,j)}\cong \widetilde{C_i}\times\widetilde{D_j}$.}

{\it Proof of claim.} We know
$$
A\times B_{\Psi(i,j)}\cong C_i\times D_j\cong W\times X\times\widetilde{C_i}\times\widetilde{D_j}\cong A\times\widetilde{C_i}\times\widetilde{D_j},
$$
\noindent so $B_{\Psi(i,j)}\cong\widetilde{C_i}\times\widetilde{D_j}$.  (See \cite[(4.3)]{LovFG}.)$\qed$

Thus $B\cong Y\times Z$.  By definition, $C\cong W\times Y$ and $D\cong X\times Z$.  $\qed$

A ``strong'' or ``strict'' version of the above theorem would be useful ({\sl a la} \cite[Proposition 3.1]{McKJC}).

\begin{lemma} Let $A$ be a finite, connected, directly irreducible poset.  Let $B$ and $P$ be posets such that $P\ne\emptyset$.  If $A\cong\mathcal C(B^P)$, then $P$ is finite and connected and $B$ is finite, connected, and directly irreducible.
\end{lemma}

\proof We have $B\ne\emptyset$.  Also $|B|\ge2$.  Note that $P$ has finitely many components.  If $P=C+D$ for posets $C,D\ne\emptyset$, then, by Proposition 1(2), $A\cong\mathcal C(B^C)\times\mathcal C(B^D)$. Therefore $A\cong\mathcal C(B^C)$, without loss of generality.  Hence $P$ is connected.

If $B=E+F$ for posets $E,F\ne\emptyset$, then by Proposition 1(3), $A\cong\mathcal C(E^P)+\mathcal C(F^P)$, a contradiction.  Hence $B$ is connected.  Thus there exist $x,y\in B$ such that $x<y$, so if $P$ is infinite, then $|P|\le|{\bf 2}^P|\le|\mathcal C(B^P)|$, a contradiction.  Hence $P$ is finite.  Also, $B$ is finite.  If $B\cong G\times H$ for posets $G$ and $H$, then by Proposition 1(4), $A\cong\mathcal C(G^P)\times\mathcal C(H^P)$, therefore without loss of generality $|\mathcal C(G^P)|=1$, so $|G|=1$. $\qed$

\begin{lemma} Let $A$ be a connected, finite, directly irreducible poset.  Let $B$ be a finite, non-empty, connected poset.  Then:
\begin{itemize}
\item[(1)] $\mathcal C(A^B)$ is connected and directly irreducible.
\item[(2)] Let $C$ and $D$ be non-empty posets.  If $\mathcal C(C^D)\cong\mathcal C(A^B)$, then $C$ is directly irreducible, finite, and connected, and $D$ is finite and connected.
\item[(3)] There exist unique (up to isomorphism) non-empty posets $E$ and $F$ with the following two properties: (a) $\mathcal C(A^B)\cong\mathcal C(E^F)$ and (b) whenever $C$ and $D$ are non-empty posets such that $\mathcal C(A^B)\cong\mathcal C(C^D)$, then there exists a poset $G$ such that $\mathcal C(E^G)\cong C$ and $G\times D\cong F$.  We can choose as $E$ any absolutely $\mathcal C$-indecomposable poset $E$ such that $\mathcal C(E^J)\cong\mathcal C(A^B)$ for some non-empty poset $J$, and we can choose that $J$ as our ``$F$.''  If $H$ is a finite, non-empty, connected poset, then $E$ and $F\times H$ are the posets that work for $\mathcal C(A^{B\times H})$.
\end{itemize}
\end{lemma}

\proof (1) The poset $\mathcal C(A^B)$ is non-trivial since $A$ is non-trivial.  By Theorem 2(3), $\mathcal C(A^B)$ is directly irreducible.  It is connected since $A$ and $B$ are connected.

(2) See the previous lemma.

(3) Using (1), take the poset $E$ given by \cite[Theorem 5.1]{McKJJ}.  There exists a non-empty finite poset $F$ such that $\mathcal C(E^F)\cong\mathcal C(A^B)$.  By (2), $F$ is connected.

If $\mathcal C(A^B)\cong\mathcal C(C^D)$, where $C$ and $D$ are non-empty posets, then by (2) $C$ and $D$ are finite and connected and $\mathcal C(E^F)\cong\mathcal C(C^D)$, so by Theorem 4, there exist non-empty, finite, connected posets $U$, $G$, $W$, and $X$ such that
$$
C\cong\mathcal C(U^G)\text{,}\ E\cong\mathcal C(U^W)\text{,}\ D\cong W\times X\text{, }\ F\cong G\times X.
$$
\noindent By absolute $\mathcal C$-indecomposability, $|W|=1$ and $E\cong U$, so the result follows.

If $E'$ and $F'$ have the same property as $E$ and $F$, then there exists a poset $G'$ such that $\mathcal C(E'^{G'})\cong E$ and $G'\times F\cong F'$.  By absolute $\mathcal C$-indecomposability, $E'\cong E$ and $|G'|=1$, so $F\cong F'$. $\qed$

\begin{theorem} Let $A$, $B$, $C$, and $D$ be finite non-empty posets such that $A$ is connected, $B$ and $D$ are non-trivial and connected, and $\mathcal C(B^A)\cong\mathcal C(D^C)$.  Then there exist finite, non-empty posets $W$, $X$, $Y$, and $Z$ such that $Z$ is non-trivial and $A\cong W\times X$, $B\cong\mathcal C(Z^Y)$, $C\cong W\times Y$, and $D\cong\mathcal C(Z^X)$.
\end{theorem}

\proof Let $B=\prod_{h\in H} B_h$ be a representation of $B$ as a product of connected, directly irreducible posets.  Let $D=\prod_{j\in J}D_j$ be a similar product.  Let $C=\sum_{i\in I} C_i$ be a decomposition into connected components.  Let $A_\ell$ ($\ell\in L$) be the finitely many pairwise non-isomorphic connected, directly irreducible posets that could arise in the factorizations of $A$ and the $C_i$ ($i\in I$).

Say $A\cong\prod_{\ell\in L} A_\ell^{k_\ell}$ where $k_\ell\ge0$ ($\ell\in L$).

For each $\ell\in L$, let $c_\ell\in\mathbb N_0$ be the maximum power of $A_\ell$ such that $A_\ell^{c_\ell}$ is a factor of all of the $C_i$ ($i\in I$).  That means, of course, that it is the highest power of $A_\ell$ in {\sl one} of the $C_i$.

Let $W=\prod_{\ell\in L} A_\ell^{\min\{c_\ell,k_\ell\}}$.  Let $X=\prod_{\ell\in L} A_\ell^{k_\ell-\min\{c_\ell,k_\ell\}}$.  Then
$$
W\times X\cong\prod_{\ell\in L} A_\ell^{k_\ell}\cong A.
$$

\quad {\it Claim 1. For all $i\in I$, $W$ is a factor of $C_i$.} 

{\it Proof of claim.} We made sure that $A_\ell^{c_\ell}$ is a factor of every $C_i$, and $W$ is a product of all of those powers of $A_\ell$ or even smaller powers. Since the different $A_\ell$ are pairwise non-isomorphic and we are working with connected posets, by Hashimoto's Refinement Theorem we have ``unique factorization,'' so $A_\ell^{c_\ell}\times A_{\ell'}^{c_{\ell'}}$ is a factor if and only if  $A_\ell^{c_\ell}$ is a factor and $A_{\ell'}^{c_{\ell'}}$ is a factor, when $\ell\ne{\ell'}$. $\qed$

By Claim 1, we may let $\widetilde{C_i}$ be $C_i$ with $W$ factored out ($i\in I$).  Let $Y=\sum_{i\in I}\widetilde{C_i}$.  Then $W\times Y\cong\sum_{i\in I} W\times\widetilde{C_i}\cong\sum_{i\in I} C_i=C$.

Note that by Proposition 1(1),
$$
\mathcal C(B^A)\cong\mathcal C\big(\mathcal C(B^X)^W\big)\ \text {and}\ \mathcal C(D^C)\cong\mathcal C\big(\mathcal C(D^Y)^W\big)
$$
\noindent and by Theorem 2(2), since $W$, being a factor of a finite, non-empty, connected poset, is finite, non-empty, and connected, $\mathcal C(B^X)\cong\mathcal C(D^Y)$.  Note that $X$ and $Y$ have no non-trivial factor in common.  (Assume for a contradiction that $X$ and $Y$ do have a non-trivial factor in common.  Then since $X$ is a factor of $A$, we may assume the common factor is $A_\ell$ for some $\ell\in L$.  As $A$ is connected, for $A_\ell$ to be a factor of $Y$, it must be a factor of $\widetilde{C_i}$ for all $i\in I$.  But then we would have pulled it out with $W$, as it were.  To be precise, since $A_\ell$ is a factor of $X$, we must have $c_\ell<k_\ell$, and so there is an $i\in I$ such that $\widetilde{C_i}$ has no factor of $A_\ell$, a contradiction.)

Write $\mathcal C(B^X)\cong\prod_{h\in H}\mathcal C(B_h^X)$ and $\mathcal C(D^Y)\cong\prod_{i\in I}\prod_{j\in J}\mathcal C(D_j^{\widetilde{C_i}}).$  By Theorem 2(3) or Lemma 7(1), $\mathcal C(B_h^X)$ ($h\in H$) and $\mathcal C(D_j^{\widetilde{C_i}})$ ($i\in I$, $j\in J$) are directly irreducible.  As they are connected, by \cite[Corollary 2]{DufHDb} (cf. \cite[Theorem 6.4]{McKGA}) there is a bijection $\Psi:I\times J\to H$ such that $\mathcal C(D_j^{\widetilde{C_i}})\cong\mathcal C(B_h^X)$ for all $(i,j)\in I\times J$.  By Theorem 4 there are finite, non-empty, connected posets $U_h$, $R_h$, $S_h$, and $T_h$ (where $h=\Psi(i,j)$) such that
$$
B_h\cong\mathcal C(U_h^{R_h}),\ D_j\cong\mathcal C(U_h^{S_h}),\ X\cong S_h\times T_h,\ \text{and}\ \widetilde{C_i}\cong R_h\times T_h.
$$
\noindent Moreover, by Theorem 2(3) or Lemma 6, $U_h$ is directly irreducible.

\quad {\it Claim 2. For each $j\in J$, let $E_j$ and $F_j$ be the posets corresponding to $D_j$ given by the part (3) of the previous lemma.  Then $X$ is a factor of $F_j$.}

{\it Proof of claim.} By part (2) of the previous lemma, $F_j$ is finite and connected.  Assume for a contradiction that $X$ does not divide $F_j$.  Then there exists $\ell\in L$ such that $A_\ell^{k_\ell-c_\ell}$ does not divide $F_j$, and so $k_\ell>c_\ell$.  Pick $i\in I$ such that $c_\ell$ is the highest power of $A_\ell$ dividing $C_i$.  Consider
$$
\mathcal C(D_j^{C_i})\cong\mathcal C(D_j^{\widetilde{C_i}\times W})\cong\mathcal C\big(\mathcal C(D_j^{\widetilde{C_i}})^W\big)\cong\mathcal C\big(\mathcal C(B_{\Psi(i,j)}^X)^W\big)\cong\mathcal C(B_{\Psi(i,j)}^{X\times W})\cong\mathcal C(B_{\Psi(i,j)}^A).
$$
\noindent By part (3) of the previous lemma, $A$ is a factor of $F_j\times C_i$ and $A_\ell^{k_\ell}$ divides $A$, so $A_\ell^{k_\ell-c_\ell}$ divides $F_j$, a contradiction. $\qed$

By Claim 2, we may let $\widetilde{F_j}$ be $F_j$ with $X$ factored out ($j\in J$).  

\quad {\it Claim 3. For all $(i,j)\in I\times J$, $B_{\Psi(i,j)}\cong\mathcal C(E_j^{\widetilde{C_i}\times\widetilde{F_j}})$.}

{\it Proof of claim.} Note that
$$
\mathcal C(B_{\Psi(i,j)}^A)\cong\mathcal C(D_j^{C_i})\cong\mathcal C\big(\mathcal C(E_j^{F_j})^{C_i}\big)\cong\mathcal C\big(\mathcal C(E_j^{X\times\widetilde{F_j}})^{W\times\widetilde{C_i}}\big)\cong\mathcal C\big(\mathcal C(E_j^{\widetilde{C_i}\times\widetilde{F_j}})^A\big)
$$
so $B_{\Psi(i,j)}\cong\mathcal C(E_j^{\widetilde{C_i}\times\widetilde{F_j}})$ by Theorem 2(2) since $A$ is connected. $\qed$

Thus
$$
B\cong\prod_{\substack{i\in I\\j\in J}}\mathcal C\big(\mathcal C(E_j^{\widetilde{F_j}})^{\widetilde{C_i}}\big)\cong\prod_{j\in J}\mathcal C\big(\mathcal C(E_j^{\widetilde{F_j}})^Y\big)\cong\mathcal C\big([\prod_{j\in J}\mathcal C(E_j^{\widetilde{F_j}})]^Y\big).
$$
\noindent Let $Z=\prod_{j\in J}\mathcal C(E_j^{\widetilde{F_j}})$.

We know
$$
D\cong\prod_{j\in J} D_j\cong\prod_{j\in J}\mathcal C(E_j^{F_j})\cong\prod_{j\in J}\mathcal C\big(\mathcal C(E_j^{\widetilde{F_j}})^X\big)\cong\mathcal C\big([\prod_{j\in J}\mathcal C(E_j^{\widetilde{F_j}})]^X\big)\cong\mathcal C(Z^X). \qed
$$

We also need a trivial extension of \cite[Lemma 2.3]{McKJJ} (although the added trivialities take up perhaps more space than is warranted).

\begin{lemma} Let $A$ and $B$ be finite posets.  Then:
\begin{itemize}
\item[(1)] $A^B=\emptyset$ if and only if $A=\emptyset$ and $B\ne\emptyset$;
\item[(2)] $A^B$ is an antichain if and only if $A$ is an antichain or $B=\emptyset$, in which case
$$
|A^B|=
\begin{cases}
|A|^c\ \text{if}\ A\ne\emptyset\ \text{or}\ B\ne\emptyset,\ \text{where $c$ is the number of connected components of}\ B\\
1\ \text{if}\ A=B=\emptyset\text{;}
\end{cases}
$$
\item[(3)] let $f,g\in A^B$ be such that $f\le g$. Assume $A\ne\emptyset$. Then $\h_{A^B}([f,g])=\h(A^B)$ if and only if $f,g\in\mathcal D(A^B)$ and $\h_A\big([f(b),g(b)]\big)=\h(A)$ for all $b\in B$;
\item[(4)] when $A\ne\emptyset$, $\h(A^B)=\h(A)|B|$.
\end{itemize}
\end{lemma}

\proof (1) This is clear.

(2) Assume $A^B$ is an antichain.  If $A$ is not an antichain and $B\ne\emptyset$, then $\{\langle a\rangle\mid a\in A\}\cong A$ is a subposet of $A^B$, a contradiction.  Conversely, if $B=\emptyset$, then $c=0$ and $|A^B|=1$, so $A^B$ is an antichain.  If $B\ne\emptyset$ but $A$ is an antichain, then for all $f,g\in A^B$, if $f<g$, there exists $b\in B$ such that $f(b)<g(b)$ in $A$, a contradiction.  Hence $A^B$ is an antichain.  For all $f\in A^B$ and for every connected component $C$ of $B$, $|f[C]|=1$, so $|A^B|=|A|^c$.

(3) This is \cite[Lemma 2.3(3)]{McKJJ} if $f<g$.  Now assume $f=g$.  If $\h_{A^B}([f,g])=\h(A^B)$, then $A^B$ is a non-empty antichain, so by (2) either $B=\emptyset$ and the forward implication is vacuous or $B\ne\emptyset$ and $A$ is a non-empty antichain, so by the above, $f,g\in\mathcal D(A^B)$ and for all $b\in B$, $f(b)=g(b)$, so $\h_A\big([f(b),g(b)]\big)=0=\h(A)$.

Conversely, suppose $f,g\in\mathcal D(A^B)$ and $0=h_A\big([f(b),g(b)]\big)=h(A)$ for all $b\in B$.  If $B=\emptyset$ then $\h(A^B)=\h({\bf 1})=0=h_{A^B}([f,g])$.  If $B\ne\emptyset$, let $b_0\in B$.  The fact $\h_A\big([f(b_0),g(b_0)]\big)=0=\h(A)$ means $A$ is an antichain.

By (1) and (2), $A^B$ is a non-empty antichain, so $\h(A^B)=0=\h_{A^B}([f,g])$.

(4) This is true by \cite[Corollary 2.2]{DufRivGH} if $A,B\ne\emptyset$.  If $A\ne\emptyset$ but $B=\emptyset$, then $\h(A^B)=\h({\bf 1})=0=\h(A)\cdot0$.$\qed$

\begin{theorem} Let $P$ and $Q$ be finite, non-empty posets such that $P$ is connected.  Assume $P^P\cong Q^Q$.   Then $Q$ is connected, and therefore $P\cong Q$.
\end{theorem}

\proof If $Q$ is connected, then $P\cong Q$ by \cite[Theorem]{DufWilGI}.  Assume now that $Q$ is disconnected. Say $Q=Q_0+D$, where $D\ne\emptyset$ and $Q_0$ is any component of $Q$ of maximum height.

We will show that $Q_0$ is the unique component of $Q$ of height $\h(Q)$, and we will show that $Q_0$ is a proper direct factor of $P$.

\quad {\it Claim 1. $P^P$ is not an antichain. Hence, $P$ is not an antichain and $|Q_0|\ne1$.}

{\it Proof of claim.} If $Q^Q\cong P^P$ is an antichain, then, by Lemma 9(2), since $P\ne\emptyset$, $P$ is an antichain. As $P$ is connected that means $|P|=1$ and $|P^P|=1$.  But by Lemma 9(2), $Q$ is an antichain and $|Q^Q|=|Q|^{|Q|}\ge2^2=4\ne1$, a contradiction.

If $P$ were an antichain, then $P^P$ would be an antichain by part (2) of the previous lemma.  If $|Q_0|=1$, then $\h(Q)=0$ and thus $Q$ is an antichain, so $Q^Q\cong P^P$ is an antichain.$\qed$

Pick $q_0\in Q_0$ of maximum height in $Q$.  Using part (3) of Lemma 9, $\langle q_0\rangle$ corresponds to $\langle p_0\rangle$ for some $p_0$ of maximum height in $P$.  Since $P$ is connected, $\mathcal C(P^P)=\{f\in P^P\mid f\equiv\langle p_0\rangle\}$.

\quad {\it Claim 2. $\mathcal C(Q_0^Q)=\{g\in Q^Q\mid g\equiv\langle q_0\rangle\}$.}

{\it Proof of claim.} Let $g\in Q^Q$ be such that $g\equiv\langle q_0\rangle$.  Then $g[Q]\subseteq Q_0$ and $g\in\mathcal C(Q_0^Q)$.

Let $h\in\mathcal C(Q_0^Q)$.  Then $h\in Q^Q$, and, in $Q_0^Q$, $h\equiv k$ for some $k\in\mathcal D(Q_0^Q)$, and $k\equiv\langle q\rangle$ for some $q\in Q_0$.  But $Q_0$ is connected, so $\langle q\rangle\equiv\langle q_0\rangle$.  Hence $h\equiv\langle q_0\rangle$.$\qed$

We conclude that $\mathcal C(P^P)\cong\mathcal C(Q_0^Q)$ via the original isomorphism. This is because we have described both $\mathcal C(P^P)$ and $\mathcal C(Q_0^Q)$ in terms just involving $P^P$ and $Q^Q$, respectively, and in the same way, up to the isomorphism, since the original isomorphism maps $\langle p_0\rangle$ to $\langle q_0\rangle$.  This proves that $Q_0$ is the unique component of $Q$ of height $\h(Q)$.  By Proposition 1(2), $\mathcal C(Q_0^Q)\cong\mathcal C(Q_0^{Q_0})\times\mathcal C(Q_0^D)\cong\mathcal C(P^P)$.  By Theorem 2(3), there exist finite, non-empty, connected posets $A_1$ and $A_2$ such that $P\cong A_1\times A_2$ and $\mathcal C(Q_0^{Q_0})\cong\mathcal C(A_1^P)$.  By Theorem 4, there exist non-empty, finite, connected posets $E$, $X$, $Y$, and $Z$ such that
$$
A_1\cong\mathcal C(E^X)\text{,}\ Q_0\cong\mathcal C(E^Y)\text{,}\ Q_0\cong X\times Z\text{,}\ P\cong Y\times Z.
$$

Since $E$ and $Y$ are connected, there exists a finite, non-empty, connected poset $F$ such that $X\cong\mathcal C(F^Y)$ by Theorem 2(3) applied to the second and third isomorphisms above.  Now suppose $|F|>1$.  Since $F$ is connected, it contains $\bf 2$ and thus $\mathcal C(F^Y)$ and hence $X$ contain ${\bf 2}^Y$.  Since the dual of $Y$, $Y^\partial$, can be embedded in ${\bf 2}^Y$, then $Y^\partial\times Z$ can be embedded in ${\bf 2}^Y\times Z$ and thus in $X\times Z\cong Q_0$.  Hence $$
\h(P)=\h(Y\times Z)=\h(Y)+\h(Z)=\h(Y^\partial)+\h(Z)=\h(Y^\partial\times Z)\le\h(Q_0)=\h(Q)
$$
\noindent and
$$
|P|=|Y\times Z|=|Y||Z|=|Y^\partial||Z|=|Y^\partial\times Z|\le|Q_0|<|Q|\text{,}
$$
\noindent so by Lemma 9(4) and Claim 1, $\h(P)\ne0$, so $\h(P^P)=|P|\h(P)<|Q|\h(Q)=\h(Q^Q)$, a contradiction.

Thus $|F|=1$, and $|X|=1$, and so $P\cong Y\times Q_0$.  If also $|Y|=1$, then $\h(P)=\h(Q)$ but $|Q|>|P|$, so $\h(Q^Q)>\h(P^P)$, a contradiction stemming from Lemma 9(4), unless $\h(P)=0$, which would make $P$ an antichain, contradicting Claim 1.  Thus $Q_0$ is a proper direct factor of $P$.

We have that $P$ and $Q_0$ are non-empty and connected, so Theorem 8 applies to the isomorphism we have already established, $\mathcal C(P^P)\cong\mathcal C(Q_0^Q)$:  There exist finite, non-empty posets $E'$, $X'$, $Y'$, and $Z'$ such that
$$
P\cong\mathcal C({E'}^{X'})\text{, }\ Q_0\cong\mathcal C({E'}^{Y'})\text{, }\ P\cong Y'\times Z'\text{, and}\ Q\cong X'\times Z'.
$$
\noindent Since $P$ is connected and non-empty, so are $Y'$ and $Z'$, and hence $Z'$ divides every component of $Q$, in particular, $Q_0$---say, $Z'\times T\cong Q_0$ for some connected, non-empty poset $T$. Again, $Y'$ is connected and so is $E'$, by Proposition 1(3), since $Q_0$ is connected. By Theorem 2(3), since $T$ is a factor of $Q_0\cong\mathcal C({E'}^{Y'})$, there exists a non-empty poset $H$ such that $T\cong\mathcal C(H^{Y'})$; $H$ is connected by Proposition 1(3).

\quad {\it Case 1. $|H|>1$}

Then since $H$ is connected, it contains $\bf 2$, and hence $T$ contains ${\bf 2}^{Y'}$, which contains $Y'^\partial$, and $Q_0\cong Z'\times T$ contains $Z'\times Y'^\partial$ so $|Q_0|\ge|Z'\times Y'^\partial|=|Z'||Y'^\partial|=|Z'||Y'|=|Z'\times Y'|=|P|$.  As we already know that $P$ properly contains $Q_0$, we have a contradiction. 

\quad {\it Case 2. $|H|=1$}

Then $|T|=1$ and $Z'\cong Q_0$, so $Q_0$ is a factor of $Q$.  But then $\h(Q_0)=\h(Q)=\h(Q_0)+\h(X')$, so $\h(X')=0$ and $X'$ is an antichain.  But if $|X'|>1$, then $Q$ has {\sl two} components of maximum height, a contradiction. $\qed$

\bigskip

Perhaps an argument like the one in \cite{FarBJb} could help us prove that $P^P\cong Q^Q$ implies $P\cong Q$ if $P$ and $Q$ are finite and non-empty and $P$ is directly irreducible.


\end{document}